# A second order positivity preserving well-balanced finite volume scheme for Euler equations with gravity for arbitrary hydrostatic equilibria

Andrea Thomann*†, Markus Zenk‡, Christian Klingenberg‡

April 26, 2018

## Abstract

We present a well-balanced finite volume solver for the compressible Euler equations with gravity where the approximate Riemann solver is derived using a relaxation approach. Besides the well-balanced property, the scheme is robust with respect to the physical admissible states. Another feature of the method is that it can maintain general stationary solutions of the hydrostatic equilibrium up to machine precision. For the first order scheme we present a well-balanced and positivity preserving second order extension using a modified minmod slope limiter. To maintain the well-balanced property, we reconstruct in equilibrium variables. Numerical examples are performed to demonstrate the accuracy, well-balanced and positivity preserving property of the presented scheme for up to 3 space dimensions.

## 1   Introduction

We consider the system of compressible Euler equations with gravity in $d$-space dimensions which is given by the following set of equations

$$\begin{aligned}\partial_t \rho + \nabla_{\mathbf{x}} \cdot (\rho \mathbf{u}) &= 0, \\ \partial_t \rho \mathbf{u} + \nabla_{\mathbf{x}} \cdot (\rho \mathbf{u} \otimes \mathbf{u} + Ip) &= -\rho \nabla_{\mathbf{x}} \Phi, \\ \partial_t E + \nabla_{\mathbf{x}} \cdot ((E+p)\mathbf{u}) &= -\rho \langle \mathbf{u}, \nabla_{\mathbf{x}} \Phi \rangle.\end{aligned} \quad (1.1)$$

---

*University of Insubria, Via Valleggio 11, 22100 Como (CO), Italy
†Marie Sklodowska Curie fellow of the Istituto Nazionale di Alta Matematica
‡Universität Würzburg, Emil-Fischer-Str. 40, 97074 Würzburg, Germany

Here, $\rho > 0$ denotes the density, $\mathbf{u} \in \mathbb{R}^d$ the velocity vector, $p$ the pressure and $E = \rho e + \frac{1}{2}\rho \|\mathbf{u}\|^2$ the total energy, where $e > 0$ is the internal energy. The function $\Phi : \mathbb{R}^d \to \mathbb{R}$ denotes a given gravitational potential. The pressure is described by a general pressure law which depends on the internal energy and specific volume $\tau = \frac{1}{\rho}$. We require for the solution $w = (\rho, \rho\mathbf{u}, E)$ the density and internal energy to be positive. That means the state vector $w$ must belong to the set $\Omega = \left\{ w \in \mathbb{R}^{d+2} \mid \rho > 0,\ e > 0 \right\}$.

The compressible Euler equations with a gravitational source term are quite important in many applications, be it in atmospheric modelling or in astrophysical stellar evolution. Typically in these applications is, that the solutions are not far from a stationary solution with zero velocity. Such special solutions are called hydrostatic equilibria. Therefore it is necessary to have a numerical method that captures those hydrostatic equilibria to machine precision, in order to resolve the evolutions near the equilibrium even when given a coarse mesh. Such numerical methods are called well-balanced methods.

In literature an exhaustive list of papers concerning well-balanced schemes can be found. We will shortly mention a few approaches. A wide range of schemes are based on finite volume methods, e.g.[1,2], central schemes e.g.[3] or discontinuous Galerkin approaches, e.g.[4–6] Challenging is also how the source term is treated. It can be included into the flux function, see e.g.[7,8] or left outside, see e.g.[6]

All those mentioned papers have in common that they can well-balance a certain class of equilibria, for example hydrostatic equilibria with constant enthalpy, see[1,2] or isothermal and polytropic equilibria, see.[7] Our paper presents a method where a given arbitrary hydrostatic equilibrium is well-balanced by a novel approach. The hydrostatic equilibrium that is chosen to well-balance is not restricted to a certain class. To approximate correctly the physically nature of the solution it is crucial that the positivity of density and total Energy is preserved which our presented scheme provably provides.

The paper is organized as follows. Section 2 is devoted to the description of hydrostatic equilibria. We will utilize the time-independent nature of stationary solutions to rewrite the derivative of the potential in terms of time-independent functions that describe the steady-state solution. In Section 3, we describe the Relaxation model that is used to derive an approximate Riemann solver. We use the results from Section 2 to achieve the well-balanced property. Section 4 is devoted to the numerical scheme which is described in a higher dimensional setup. A Godunov-type first order scheme is given for which a second order extension is described. It is based on linear reconstruction on cell averages and, to keep the well-balanced property, on a transformation into equilibrium variables.

For the resulting scheme, the main properties, which are the robustness and the well-balancednes of the scheme, are proven in Section 5. It is followed by a section with numerical results to validate the main properties given in Section 5. We start with the well-balancedness by calculating isothermal, polytropic and a general stationary equilibrium and give the $L^1$-error. To show the accuracy of the second order scheme, we consider two analytical solutions of the Euler equations with gravity. One is already known from literature[6] and a novel three dimensional one derived from an isothermal approach. Additionally, the convergence of the source term is shown by balancing a general stationary state against an isothermal equilibrium. We also perform perturbation test cases like a Rayleigh-Taylor instability[9] and a perturbation in pressure upon a hydrostatic equilibrium.

To demonstrate the positivity preserving property, we perform a double rarefraction test based on an isothermal setting similar to the Einfeldt rarefraction test.[10]

A section of conclusion completes this paper.

## 2 Hydrostatic states

In the following, we will focus on steady states at rest, which are solutions of

$$\mathbf{u} = 0, \tag{2.1}$$

$$\nabla_\mathbf{x} p = -\rho \nabla_\mathbf{x} \Phi. \tag{2.2}$$

Let a hydrostatic stationary solution be given through $\bar{u}, \bar{\rho}$ and $\bar{p}$. Then $\bar{\rho}$ and $\bar{p}$ are time-independent. Following,[11] we can write the hydrostatic solution as just space dependent functions

$$\bar{\rho} = \alpha(\mathbf{x}) \text{ and } \bar{p} = \beta(\mathbf{x}). \tag{2.3}$$

Since the density and the pressure are strictly positive, we also require $\alpha, \beta > 0$. These functions must satisfy the hydrostatic equation (2.2) which leads to

$$\nabla_\mathbf{x} \beta(\mathbf{x}) = -\alpha(\mathbf{x}) \nabla_\mathbf{x} \Phi(\mathbf{x}). \tag{2.4}$$

Using this relation, we can find the following expression for the gradient of the gravitational potential

$$\nabla_\mathbf{x} \Phi(\mathbf{x}) = -\frac{\nabla_\mathbf{x} \beta(\mathbf{x})}{\alpha(\mathbf{x})}. \tag{2.5}$$

In the following, we will consider the hydrostatic equation with the rewritten potential gradient

$$\nabla_\mathbf{x} p = \rho \frac{\nabla_\mathbf{x} \beta}{\alpha}. \tag{2.6}$$

For illustration we give some examples how the functions $\alpha$ and $\beta$ can be found for an isothermal and polytropic equilibrium.

### 2.1 Isothermal atmosphere

Using the ideal gas law $p = \rho RT$ and the isothermal condition $T = T_c = \text{const.}$, equation (2.2) has an analytical solution given by

$$\bar{p}(\mathbf{x}) = p_c \exp\left(-\frac{\Phi(\mathbf{x})}{RT_c}\right), \tag{2.7}$$

with $p_c = \rho_c R T_c$. The functions $\alpha, \beta$ are then given by

$$\alpha(\mathbf{x}) = \rho_c \exp\left(-\frac{\Phi(\mathbf{x})}{RT_c}\right), \tag{2.8}$$

$$\beta(\mathbf{x}) = p_c \exp\left(-\frac{\Phi(\mathbf{x})}{RT_c}\right). \tag{2.9}$$

## 2.2 Polytropic atmosphere

The polytropic atmosphere is characterized by the relation $\bar{p} = s\bar{\rho}^\nu$, where $s, \nu$ are positive constants and $\nu > 1$. Inserting the polytropic relation into (2.2) and differentiation gives

$$s\nu\bar{\rho}^{\nu-2}\nabla_\mathbf{x}\bar{\rho} = \nabla_\mathbf{x}\Phi. \tag{2.10}$$

Integrating this equation and rewriting leads to

$$\bar{\rho}(\mathbf{x}) = \rho_c\left(1 + \frac{\nu-1}{s\nu\rho_c^{\nu-1}}(\Phi_c - \Phi(\mathbf{x}))\right)^{\frac{1}{\nu-1}}. \tag{2.11}$$

Thus the functions $\alpha, \beta$ are given by

$$\alpha(\mathbf{x}) = \rho_c\left[1 + \frac{\nu-1}{s\nu\rho_c^{\nu-1}}(\Phi_c - \Phi(\mathbf{x}))\right]^{1/(\nu-1)}, \tag{2.12}$$

$$\beta(\mathbf{x}) = s\,\alpha(\mathbf{x})^\nu. \tag{2.13}$$

# 3 Relaxation model

We consider the following relaxation model, for simplicity in one spatial direction, as derived in.[7]

$$\begin{aligned}
\partial_t \rho + \partial_{x_1} \rho u_1 &= 0, \\
\partial_t \rho u_1 + \partial_{x_1}(\rho u_1^2 + \pi) &= -\rho\,\partial_{x_1} Z, \\
\partial_t \rho u_2 + \partial_{x_1}(\rho u_1 u_2) &= 0, \\
\partial_t \rho u_3 + \partial_{x_1}(\rho u_1 u_3) &= 0, \\
\partial_t E + \partial_{x_1}(E + \pi)u_1 &= -\rho\,u_1\partial_{x_1} Z, \\
\partial_t \rho\pi + \partial_{x_1}(\rho\pi + a^2)u_1 &= \frac{\rho}{\epsilon}(p(\tau, e) - \pi), \\
\partial_t \rho Z + \partial_{x_1}\rho Z u_1 &= \frac{\rho}{\epsilon}(\Phi - Z).
\end{aligned} \tag{3.1}$$

The pressure is approximated by a new variable $\pi$ following the Suliciu relaxation approach described in,[12] where $a > 0$ denotes the relaxation parameter. Additionally, the potential $\Phi$ is approximated by a new variable $Z$ which is transported with velocity $u_1$.

Analogously to (1.1), we define the state vector $W = (\rho, \rho\mathbf{u}, E, \rho\pi, \rho Z)$ which belongs to $\Omega_W = \{W \in \mathbb{R}^{4+d}, \rho > 0, e > 0\}$. For a given gravity function $\Phi$, an relaxation equilibrium state for model (3.1) is defined by

$$W^{eq} = (\rho, \rho\mathbf{u}, E, \rho p(\tau, e), \rho\Phi)^T. \tag{3.2}$$

In the following, we sum up the properties of the model (3.1). A more detailed description can be found in.[7,12] The eigenvalues of the system in primitive variables are $\lambda^\pm = u_1 \pm \frac{a}{\rho}$ and $\lambda^{\mathbf{u}} = u_1$, where the eigenvalue $\lambda^{\mathbf{u}}$ has multiplicity five. One finds the fields associated to the eigenvalues are linearly degenerate and the Riemann invariants with respect to $\lambda^\pm$ are

$$I_1^\pm = u_1 \pm \frac{a}{\rho}, \ I_2^\pm = \pi \mp au, \ I_3^\pm = e - \frac{\pi^2}{2a^2}, \ I_4^\pm = Z, \ I_{5,6}^\pm = u_{2,3} \qquad (3.3)$$

and with respect to $\lambda^{\mathbf{u}}$

$$I^{\mathbf{u}} = u_1. \qquad (3.4)$$

Let us consider as initial data a Riemann problem at $x_1 = 0$ with two constant values $W_L, W_R$ which are separated by a discontinuity at $x_1 = 0$

$$W_0(x_1) = \begin{cases} W_L & x_1 < 0 \\ W_R & x_1 > 0. \end{cases} \qquad (3.5)$$

The solution consists of four constant states separated by contact discontinuities and has the following structure

$$W_R\left(\frac{x_1}{t}; W_L, W_R\right) = \begin{cases} W_L & \frac{x}{t} < \lambda^- \\ W_L^* & \lambda^- < \frac{x_1}{t} < \lambda^u \\ W_R^* & \lambda^u < \frac{x_1}{t} < \lambda^+ \\ W_R & \lambda^+ < \frac{x_1}{t} \end{cases}, \qquad (3.6)$$

where $W_L^*, W_R^*$ denote the intermediate states. The resulting Riemann problem consists of $2(4+d)$ unknowns, $4+d$ for each intermediate state $W_{L,R}^*$. To solve the Riemann problem $2(4+d)$ relations are needed but one obtains just $2(3+d)+1$ relations from the Riemann invariants (3.3) and (3.4). This leaves us with one degree of freedom to choose the missing relation such that the resulting scheme has the well-balanced property.

In relaxation variables, the hydrostatic equilibrium (2.2) in $x_1$ direction is given by

$$\partial_{x_1}\pi = -\rho\partial_{x_1}Z. \qquad (3.7)$$

Let

$$\pi_R^* - \pi_L^* = \overline{S} \qquad (3.8)$$

be a discretization of (3.7), where $\overline{S}$ is a discretization of the source term.

Using this relation in addition to the relations gained from the Riemann invariants, the intermediate states $W_{L,R}^*$ can be determined. They are given

by

$$\frac{1}{\rho_L^*} = \frac{1}{\rho_L} + \frac{1}{a}(u^* - u_L), \quad \frac{1}{\rho_R^*} = \frac{1}{\rho_R} + \frac{1}{a}(u_R - u^*), \tag{3.9}$$

$$u^* = u_L^* = u_R^* = \frac{1}{2}(u_L + u_R) - \frac{1}{2a}\left(\pi_R - \pi_L - \overline{S}\right), \tag{3.10}$$

$$\pi_L^* = \pi_L + a(u_L - u^*), \quad \pi_R^* = \pi_R + a(u^* - u_R), \tag{3.11}$$

$$e_L^* = e_L + \frac{1}{2a^2}\left({\pi_L^*}^2 - \pi_L^2\right), \quad e_R^* = e_R + \frac{1}{2a^2}\left({\pi_R^*}^2 - \pi_R^2\right), \tag{3.12}$$

$$Z_L^* = Z_L, \quad Z_R^* = Z_R. \tag{3.13}$$

For more details on the computations see.[7]

Thus, the Riemann problem of the relaxation system completed by relation (3.8) has a unique solution of the form (3.6). An alternative approach which leads to the same approximate Riemann solver can be found in.[13]

Now it remains to give an explicit form of the source term $\overline{S}$. We choose the following definition

$$\overline{S}(W_L, W_R) = \left(\frac{\rho_R}{\alpha_R}\beta_R - \frac{\rho_L}{\alpha_L}\beta_L\right) - \frac{1}{2}\left(\frac{\rho_R}{\alpha_R} - \frac{\rho_L}{\alpha_L}\right)(\beta_R + \beta_L). \tag{3.14}$$

This formula can be derived by partial integration and Taylor expansion up to $\mathcal{O}(\Delta x^3)$.

The following result shows the well-balanced property of the defined approximative Riemann solver.

**Lemma 1** *Let $W_L, W_R$ be given by*

$$\mathbf{u} = 0, \quad \frac{\rho_{L,R}}{\alpha_{L,R}} = const., \quad \frac{p_{L,R}}{\beta_{L,R}} = const.. \tag{3.15}$$

Then follows for the intermediate states $W_L^* = W_L$ and $W_R^* = W_R$. That means the approximate Riemann solver is at rest.

**Proof**. Let $w_L, w_R$ be given satisfying (3.15). If $u^* = 0$ we have from the eqs. (3.9) and (3.11) to (3.13) together with $u_{L,R} = 0$ that $W_L^* = W_L$ and $W_R^* = W_R$ and thus (3.15) is fulfilled for $w^{eq}(\frac{x}{t}; w_L, w_R)$. It remains to show, that $u^* = 0$. In equilibrium we have $\frac{\rho_R}{\alpha_R} - \frac{\rho_L}{\alpha_L} = 0$, $\frac{\rho_{L,R}}{\alpha_{L,R}} = 1$ and the given formula for $\overline{S}$ simplifies to

$$\overline{S} = \frac{\rho_R}{\alpha_R}\beta_R - \frac{\rho_L}{\alpha_L}\beta_L = \beta_R - \beta_L = \pi_R - \pi_L.$$

Since $\pi_R - \pi_L = \overline{S}$ all terms on the right side in equation (3.10) vanish and we get $u^* = 0$. □

We want to conclude this section by mentioning some additional properties of the above defined Riemann solver following.[7]

**Lemma 2** *Given initial data $w_L, w_R \in \Omega$, then for the relaxation parameter $a$ large enough there is $w^{eq}(\frac{x}{t}; w_L, w_R) \in \Omega$; i.e the approximate Riemann solver is robust with respect to the positivity of density $\rho$ and internal energy $e$.*

**Proof.** The positivity of the density $\rho^*_{L,R}$ follows directly from the order of the eigenvalues $u_L - \frac{a}{\rho_L} < u^*$ from which we get $u_L - u^* < \frac{a}{\rho_L}$. Using the Riemann invariant $u - \frac{a}{\rho}$ we have

$$\frac{1}{\rho^*_L} = \frac{1}{\rho_L} + \frac{1}{a}(u^* - u_L) > \frac{1}{\rho_L} - \frac{1}{\rho_L} = 0.$$

The positivity of $\rho^*_R$ can be shown analogously.

For the positivity of the internal energy we have

$$e^*_L = e_L + \frac{1}{8a^2}(3\pi_L + \pi_R - \overline{S})(\pi_R - \pi_L - \overline{S}) + \frac{1}{4a}(\pi_R + \pi_L - \overline{S})(u_L - u_R) + \frac{1}{8}(u_L + u_R)^2.$$

This formula contains only positive terms or terms which can be controlled if $a$ is chosen sufficiently large. For $e^*_R$ there can be found a similar formula. □

**Remark 3** *Let $\partial_t \rho F(\eta) + \partial_x F(\eta)u \leq 0$ an entropy inequality for the Euler equations with gravity, where $\eta(\tau, e)$ denotes a specific entropy. Then the approximate Riemann solver is consistent with the entropy inequality. This follows directly from Theorem 8 from,[7] since the given proof there is independent of the source term and thus it can be directly applied here.*

## 4 Numerical scheme

In this section, we describe the numerical scheme associated with the approximative Riemann solver developed above. Consider at first the general form of the Euler equations (1.1) as

$$w_t + \nabla_{\mathbf{x}} \cdot f(w) = S(w). \tag{4.1}$$

Furthermore consider the relaxation system derived in section 3 as

$$W_t + \nabla_{\mathbf{x}} \cdot F(W) = S(W) + \frac{1}{\epsilon}R(W), \tag{4.2}$$

where we can decompose the flux function of the relaxation system as

$$F = \begin{pmatrix} F_1 \\ F_2 \end{pmatrix},$$

$$F_1 = \begin{pmatrix} \rho u_1 & \ldots & \rho u_d \\ \rho u_1^2 + \pi & \ldots & \rho u_1 u_d \\ \vdots & \ddots & \vdots \\ \rho u_1 u_d & \ldots & \rho u_d^2 + \pi \\ u_1(E + \pi) & \ldots & u_d(E + \pi) \end{pmatrix}, \tag{4.3}$$

$$F_2 = \begin{pmatrix} \rho u_1 \pi + a^2 u_1 & \ldots & \rho u_d \pi + a^2 u_d \\ \rho u_1 Z & \ldots & \rho u_d Z \end{pmatrix}.$$

We would like to note that the suggested relaxation system is consistent with (1.1) in the way that $F_1(W^{eq}) = f(w)$.

For the space discretization we construct a cartesian mesh by defining a family of stepsizes $\Delta x_i$ for $i \in \{1, ..., d\}$ for each dimension and from that a family of points

$$\mathbf{x}_I = [x_{1,I_1}, ..., x_{d,I_d}] = [x_{1,0} + I_1 \Delta x_1, ..., x_{d,0} + I_d \Delta x_d], \qquad (4.4)$$

where $I \in \mathcal{I}$ is a multi-index defined as $I = [I_1, ...., I_d]$. To simplify notations we further define

$$I(i, k) = [I_1, ..., I_i + k, ..., I_d]. \qquad (4.5)$$

With this, we define the cells

$$C_I = [x_{1,I(1,-\frac{1}{2})}, x_{1,I(1,\frac{1}{2})}] \times ... \times [x_{d,I(d,-\frac{1}{2})}, x_{d,I(d,\frac{1}{2})}]. \qquad (4.6)$$

The time discretization on the interval $[0, T]$ is given by $t^{n+1} = t^n + \Delta t$ where $\Delta t > 0$ denotes the length of the time step restricted by a CFL condition. As standard in the finite volume setting, we compute approximations of the cells averages

$$w_I^n = \frac{1}{\Pi_{i=1}^d \Delta x_i} \int_{C_I} w(\mathbf{x}, t^n) d\mathbf{x} \qquad (4.7)$$

Furthermore we consider data to be in a hydrostatic equilibrium when

$$\mathbf{u}_I = \mathbf{0}, \quad \rho_I = \alpha_I, \quad p_I = \beta_I \qquad (4.8)$$

holds for all $I \in \mathcal{I}$.

## 4.1 First order scheme

As standard in a finite volume scheme, we find after integration of (4.1) over the cell $C_I$ and over the time intervall $[t^n, t^n + \Delta t]$, the evolution of the data as

$$w_I^{n+1} = w_I^n - \Delta t \sum_{i=1}^d \frac{1}{\Delta x_i} \mathbf{n}_i \cdot \left( F_{I(i,-\frac{1}{2})}^+ - F_{I(i,\frac{1}{2})}^- \right), \qquad (4.9)$$

where $F$ denotes the numerical flux and $\mathbf{n}_i$ the $i-th$ unit vector. As in a Godunov-type scheme, the fluxes are computed by considering the Riemann problem at the cell interfaces. In our case, we specifically use the relaxation system (4.2). Therefore we consider the Riemann problem as

$$\begin{aligned} W_t + F_{i,x_i} &= S_i, \\ W_0(x_i) &= \begin{cases} W_I^{eq}(w_I) & \text{if } x_i < 0 \\ W_{I(i,1)}^{eq}(w_{I(i,1)}) & \text{if } x_i > 0 \end{cases}, \end{aligned} \qquad (4.10)$$

where $F_i$ denotes the $i-th$ row of $F$ and $S_i = (0, \rho \mathbf{n}_i \odot \nabla_{\mathbf{x}} \Phi, \rho \langle \mathbf{u}, \mathbf{n}_i \odot \nabla_{\mathbf{x}} \Phi \rangle)^T$ evaluated over the interface between the cells $C_I$ and $C_{I(i,1)}$. Given (3.6), we define the numerical fluxes as follows

$$(F^-_{I(i,\frac{1}{2})}, F^+_{I(i,\frac{1}{2})}) = \begin{cases} (\mathbf{n}_i \cdot F_1(W_L), \mathbf{n}_i \cdot F_1(W_L) + S_{I(i,\frac{1}{2})}) & \text{if } \lambda^- > 0 \\ (\mathbf{n}_i \cdot F_1(W_L^*), \mathbf{n}_i \cdot F_1(W_L^*) + S_{I(i,\frac{1}{2})}) & \text{if } \lambda^{\mathbf{u}} > 0 > \lambda^- \\ (\mathbf{n}_i \cdot F_1(W_L^*), \mathbf{n}_i \cdot F_1(W_R^*)) & \text{if } \lambda^{\mathbf{u}} = 0 \\ (\mathbf{n}_i \cdot F_1(W_R^*) - S_{I(i,\frac{1}{2})}, \mathbf{n}_i \cdot F_1(W_R^*)) & \text{if } \lambda^+ > 0 > \lambda^{\mathbf{u}} \\ (\mathbf{n}_i \cdot F_1(W_R) - S_{I(i,\frac{1}{2})}, \mathbf{n}_i \cdot F_1(W_R)) & \text{if } \lambda^+ < 0 \end{cases},$$
(4.11)

where $S_{I(i,\frac{1}{2})} = (0, \mathbf{n}_i \overline{S}_{I(i,\frac{1}{2})}, \lambda^{\mathbf{u}} \overline{S}_{I(i,\frac{1}{2})})^T$ and $\overline{S}_{I(i,\frac{1}{2})} = \overline{S}(W_I, W_{I(i,1)})$.

It should be remarked, that in general $F^-_{I(i,\frac{1}{2})} \neq F^+_{I(i,\frac{1}{2})}$, since we include the source term into the flux definition. Moreover we are only advancing the physical variables and do not consider the update on the relaxation variables since we consider the initial condition to be at the relaxation equilibrium.

## 4.2 Second order in space

Concerning the spatial order, we extend the first order scheme to second order by constructing piecewise linear functions in the primitive variables $w^p = (\rho, \mathbf{u}, p)$. We make use of the structure of the cartesian mesh and reconstruct along each dimension separately. Thus, we consider a linear function in $x_i$ in the $i - th$ spatial dimension as

$$w^p_{I,i}(x_i) = w^p_I + \sigma_{I,i}(x_i - x_{i,I_i}),$$
(4.12)

where $\sigma_{I,i} \in \mathbb{R}^{d+2}$. The slopes $\sigma_{I,i}$ are computed from the neighbouring cell along that dimension, i.e.

$$\sigma_{I,i} = \text{limiter}(w^p_{I(i,-1)}, w^p_I, w^p_{I(i,1)}).$$
(4.13)

The function *limiter* can be any consistent limiter in the sense that from

$$w^p_{I(i,-1)} = w^p_I = w^p_{I(i,1)} \quad \text{it follows} \quad \sigma_{I,i} = 0,$$
(4.14)

where (4.14) is to be understood component-wise.

In order to derive a well-balanced scheme, we choose to modify this procedure by adjusting the reconstruction in pressure. We define the following transformation

$$\begin{aligned} q_{I,I(i,-1)} &= p_{I(i,-1)} + \overline{S}_{I(i,-\frac{1}{2})}, \\ q_{I,I} &= p_I, \\ q_{I,I(i,1)} &= p_{I(i,1)} - \overline{S}_{I(i,\frac{1}{2})}, \end{aligned}$$
(4.15)

and with this the set $\bar{w}^p_{I,i} = (\rho, \mathbf{u}, q)$. We would like to emphasize that this new set has to be computed for each stencil on which the reconstruction is considered.

Then the interface values are computed as

$$w^{p,\mp}_{I,i} = \bar{w}^{p,\mp}_{I,i} = w^p_I \pm \sigma_{I,i} \frac{\Delta x_i}{2}.$$
(4.16)

The initial condition for the Riemann problem at the interface $I(i, \frac{1}{2})$ is then defined as

$$W_0(x_i) = \begin{cases} W^{eq}(w_{I,i}^{p,-}) & \text{if } x_i < 0 \\ W^{eq}(w_{I(i,1),i}^{p,+}) & \text{if } x_i > 0 \end{cases} \quad (4.17)$$

and

$$\overline{S}_{I(i,\frac{1}{2})} = \overline{S}(W_I^{eq}, W_{I(i,1)}^{eq}). \quad (4.18)$$

**Remark 4** *From the reprojection in (4.16) we see that when we compute a slope for the new values $q$ we are in fact finding a slope for the pressure since $q_{I,I}$ and $p_I$ coincide on the cell $C_I$.*

**Remark 5** *For the well balance property it is crucial that the source term, as denoted in equation (4) for the source discretization, is evaluated using the cell average values.*

We would like to emphasize a crucial property of the new reconstruction in the following Lemma.

**Lemma 6** *Let the initial data be given in hydrostatic equilibrium (4.8) for all $I \in \mathcal{I}$. Then after applying the reconstruction (4.12) - (4.16) we have*

$$\mathbf{u}_{I,i}^{\mp} = \mathbf{0}, \quad p_{I,i}^{\mp} = p_I. \quad (4.19)$$

**Proof.** The reconstruction of the velocities follows from the consistency of the limiter (4.14). We now want to analyse the reconstruction in $q$. Given (4.15) and (3.14), we compute the pressure transformation

$$\begin{aligned} q_{I,I(i,-1)} &= p_{I(i,-1)} + \overline{S}_{I(i,-\frac{1}{2})} \\ &= p_{I(i,-1)} + \left( \frac{\rho_I}{\alpha_I} \beta_I - \frac{\rho_{I(i,-1)}}{\alpha_{I(i,-1)}} \beta_{I(i,-1)} \right) \\ &\quad - \frac{1}{2} \left( \frac{\rho_I}{\alpha_I} - \frac{\rho_{I(i,-1)}}{\alpha_{I(i,-1)}} \right) (\beta_I + \beta_{I(i,-1)}) \\ &= p_{I(i,-1)} + (\beta_I - \beta_{I(i,-1)}) = p_I. \end{aligned} \quad (4.20)$$

Analogously we get $q_{I,I(i,1)} = p_I$. Therefore the result follows again from the consistency of the applied limiter (4.14). □

### 4.3 Second order in time

We follow the lines of Berthon, see,[14] and consider the first order time update formula (4.9) as

$$w_I^{n+1} = w_I^n - \Delta_t H(w_\mathcal{I}^n). \quad (4.21)$$

Based on this, we construct the following two stage time integration

$$\begin{aligned} \bar{w}_I &= w_I^n - \Delta t_1 H(w_\mathcal{I}^n), \\ \bar{\bar{w}}_I &= \bar{w}_I - \Delta t_2 H(\bar{w}_\mathcal{I}), \\ w_I^{n+1} &= w_I^n (1 - \gamma) + \gamma \bar{\bar{w}}_I, \end{aligned} \quad (4.22)$$

where
$$\gamma = \frac{2\Delta t_1 \Delta t_2}{(\Delta t_1 + \Delta t_2)^2}.$$

The total time increment is given by
$$\Delta t = \frac{2\Delta t_1 \Delta t_2}{\Delta t_1 + \Delta t_2} \tag{4.23}$$

**Remark 7** *This time integrator has the advantage that the CFL criteria can be met for every single stage independently and the total time increment is a result of the possibly different increments in the single stages. Moreover it is straightforward to see that $\gamma \in (0, 1)$. This makes the final update a convex combination of the values $w_I^n$ and $\bar{\bar{w}}_I$.*

# 5 Properties of the numerical scheme

In this section we would like to state the central properties of the first and second order scheme derived above.

**Theorem 8** *Let the initial data be given in a hydrostatic equilibrium (4.8) for all $I \in \mathcal{I}$. Then the first order scheme is well-balanced.*

**Proof**. We know from Lemma 1 that the approximate Riemann solver at the cell interfaces is at rest. Therefore with the definition of the numerical fluxes (4.11) we have
$$F^{\mp}_{I(i,\pm\frac{1}{2})} = F(w_I). \tag{5.1}$$

Using (5.1) in the update formula (4.9) we get
$$w_I^{n+1} = w_I^n - \Delta t \sum_{i=1}^d \frac{1}{\Delta x_i} \mathbf{n}_i \cdot \underbrace{\left( F^+_{I(i,-\frac{1}{2})} - F^-_{I(i,\frac{1}{2})} \right)}_{=0} = w_I^n \tag{5.2}$$

and thus the first order scheme is well-balanced. □

**Theorem 9** *Let the initial data be given in a hydrostatic equilibrium (4.8) for all $I \in \mathcal{I}$. Then the second order scheme is well-balanced.*

**Proof**. First, we want to address the analysis of the interface values. From Lemma 6, we have
$$\mathbf{u}^{\mp}_{I,i} = \mathbf{0}, \quad p^{\mp}_{I,i} = p_I. \tag{5.3}$$

Furthermore by using the source average as defined in (4.18), we can proof again the results of Lemma 1 and conclude that the approximate Riemann solver stays at rest. Following the lines of the proof of Theorem 8, we conclude for the second order time update given by (4.22) we have the following implication
$$\bar{w}_I = w_I^n, \quad \bar{\bar{w}}_I = \bar{w}_I, \quad \Rightarrow w_I^{n+1} = w_I^n \tag{5.4}$$

and thus the second order scheme is well-balanced. □

**Theorem 10** *Let $w_I \in \Omega$, then under the CFL condition $\frac{\Delta t}{\Delta x_i} \max |\lambda_i| < \frac{1}{2d}$ the first order scheme is robust.*

**Proof.** We begin the proof with the analysis of the special case of one spatial dimension. Then the update (4.9) can be rewritten as

$$w_I^{n+1} = \frac{1}{\Delta x_1} \left( \int_{x_{I(1,-\frac{1}{2})}}^{x_I} W_{\mathcal{R}}^{(\rho, \rho \mathbf{u}, E)} \left( \frac{x}{t^{n+1}}, W_{I(1,-1)}, W_I \right) dx \right.$$
$$\left. + \int_{x_I}^{x_{I(1,\frac{1}{2})}} W_{\mathcal{R}}^{(\rho, \rho \mathbf{u}, E)} \left( \frac{x}{t^{n+1}}, W_I, W_{I(1,1)} \right) dx \right), \quad (5.5)$$

see also.[7] We know that $W_{\mathcal{R}}(\frac{x}{t}, W_{I(1,-1)}, W_I), W_{\mathcal{R}}(\frac{x}{t}, W_I, W_{I(1,1)}) \in \Omega$ from Lemma 2 and by the convexity of $\Omega$ and (5.5) we have $w_I^{n+1} \in \Omega$. Now we extend the analysis to higher space dimensions. For this we can rewrite the update formula (4.9) as

$$w_I^{n+1} = w_I^n - \Delta t \sum_{i=1}^d \frac{1}{\Delta x_i} \mathbf{n}_i \cdot \left( F_{I(i,-\frac{1}{2})}^+ - F_{I(i,\frac{1}{2})}^- \right)$$
$$= \frac{1}{d} \sum_{i=1}^d \underbrace{\left[ w_I^n - \frac{d \Delta t}{\Delta x_i} \mathbf{n}_i \cdot \left( F_{I(i,-\frac{1}{2})}^+ - F_{I(i,\frac{1}{2})}^- \right) \right]}_{\stackrel{(5.5)}{\in} \Omega} \quad (5.6)$$

under the CFL condition $\frac{1}{2d}$. Therefore by convexity of (5.6), we have $w_I^{n+1} \in \Omega$. □

Now we would like to tackle the issue of the robustness for the second order method. If a reconstruction in conservative variables is applied, then for the robustness it is sufficient to check the robustness of the reconstructed states at the interface, see Bouchut.[15] However we decide to reconstruct in primitive variables. In this case we can use the results of Berthon[14] in order to ensure the robustness. We give the following Theorem as a summary of the work in[14] within our context.

**Theorem 11** *Let $w_I \in \Omega$, then under the CFL condition $\frac{\Delta t}{\Delta x_i} \max |\lambda_i| < \frac{1}{3} \frac{1}{2d}$ the second order scheme is robust.*

**Proof.** The proof is a straightforward application of the results proven in section 3.1 in.[14] □

# 6 Numerical results

To illustrate the properties presented in Section 5, we present several numerical experiments. In all numerical experiments performed with the second order scheme, we use the minmod limiter defined as

$$\text{minmod}(x, y) = \begin{cases} \min(x, y) & \text{if } x, y \geq 0 \\ \max(x, y) & \text{if } x, y \leq 0 \\ 0 & \text{otherwise} \end{cases} \quad (6.1)$$

and apply the slope limiting from[14] to enusure the robustness of the scheme. All computations are carried out on double precision and errors are given in the $L^1$-norm.

## 6.1 Well-balanced tests

For the well-balanced tests, we consider different stationary solutions in one, two and three space dimensions using the second order scheme. The computations are performed on a uniform grid on the domain $[0, 1]^d$, where $d$ denotes the number of dimensions, up to a final time $T_f = 1.0$.

As a first example, we consider a isothermal hydrostatic atmosphere in three dimensions, see Section 2.1, with $\Phi(\mathbf{x}) = \frac{1}{2}(x_1^2 + x_2^2 + x_3^2)$.

Then a polytropic atmosphere in two space dimensions according to Section 2.2 with $\Phi(\mathbf{x}) = x_1 + x_2$ is considered.

Next, we consider a non-polytropic and non-isothermal stationary state with the potential $\Phi(\mathbf{x}) = -\sum_{j=1}^{d} \sin(2\pi x_j)$ on the domain $[0, 1]^d$ with periodic boundary conditions. A steady state solution for this potential is given by

$$\rho(\mathbf{x}) = c_\rho - 2\, \Phi(\mathbf{x}), \tag{6.2}$$

$$\mathbf{u} = 0, \tag{6.3}$$

$$p(\mathbf{x}) = c_p - c_\rho\, \Phi(\mathbf{x}) - \frac{1}{2} \sum_{j=1}^{d} \cos(4\pi x_j) + \sum_{i,j=1, j>i}^{d} \sin(2\pi x_i) \sin(2\pi x_j). \tag{6.4}$$

For this general stationary state calculations are performed in one, two and three space dimensions. It can be seen from Tables 2, 3 and 4 that the stationary states are preserved on machine precision.

## 6.2 Accuracy

### 6.2.1 Source term

However the equilibrium solution (6.2) can be preserved on machine precision, we want to balance it against an isothermal equilibrium using the first order scheme to demonstrate the accuracy of the source term discretization. For the given discretization $\overline{S}$ one finds with a straightforward computation using Taylor expansions that $p_R - p_L = \overline{S} + \mathcal{O}(\Delta x^3)$.

From Table 1, we see that the first order scheme with $\overline{S}$ converges with order 2. This loss of one order of accuracy can be explained that the source term is included into the flux function, see (4.11), which leads to the reduce of accuracy.

### 6.2.2 Exact solutions

To demonstrate that the second order extension of the first order scheme has the expected accuracy, we compare the numerical solution to an exact solution of the Euler equations with gravity. An exact solution of (1.1) is given by

$$\rho(\mathbf{x}, t) = \exp\left(\frac{1}{RT}\left(\frac{1}{2}\sum_{i=1}^{d} u_i^2 - \sum_{i=1}^{d}\frac{\kappa_i}{\eta_i}\cos(\eta_i t)\Phi_{x_i} - \Phi\right)\right),$$
$$u_i(t) = \kappa_i \sin(\eta_i t),$$
$$p(\mathbf{x}, t) = RT\rho(\mathbf{x}, t), \quad (6.5)$$
$$\Phi(\mathbf{x}) = \frac{1}{2}\sum_{i=1}^{d} \eta_i^2 x_i^2,$$

where the constants $\kappa_i > 0$ denote the amplitude of the velocities $u_i$ and $\eta_i > 0$ are scaling constants. The velocities are time dependent and for $\mathbf{u} = 0$, (6.5) is an isothermal hydrostatic equilibrium. We have set $\alpha$ and $\beta$ according to the isothermal equilibrium solution. It can be seen from Table 6 that the convergence rates are approaching 2.

The calculations are performed on the domain $[0,1]^3$ starting with 25 cells using exact boundary conditions. Another exact solution in two space dimensions taken from[6] is given by

$$\rho(\mathbf{x}, t) = 1 + 0.2\sin(\pi(x_1 + x_2 - t(u_{1_0} + u_{2_0}))),$$
$$u_1(t) = u_{1_0},$$
$$u_2(t) = u_{2_0},$$
$$p(\mathbf{x}, t) = p_0 + t(u_{1_0} + u_{2_0} - (x_1 + x_2)) + 0.2\cos(\pi(x_1 + x_2 - t(u_{1_0} + u_{2_0})))/\pi.$$

For the parameters, we choose $u_{1_0} = 20, u_{2_0} = 20$ and $p_0 = 4.5$. As $\alpha$ and $\beta$, we choose the density and pressure for $u_{1_0} = 0$ and $u_{2_0} = 0$ which is a stationary state of the Euler equations with gravity. The computational domain is $[0, 1]^2$ and the computations are performed with exact boundary conditions up to a final time $T_f = 0.01$. It can be seen from Table 5 that the convergence rates throughout all variables are around 1.9.

## 6.3 Evolution of small perturbations

In the following numerical test, the evolution of a small perturbation added to an initial isothermal hydrostatic solution is investigated, see.[9] The initial values on the domain $[0, 1]$ are given by

$$\Phi(x) = x,$$
$$\rho(x) = \exp(-\Phi(x)),$$
$$p(x) = \exp(-\Phi(x)) + 0.01\exp(-100(x - 0.5)^2),$$

where the pressure is perturbed by a Gauß function centered in $x = 0.5$. The solution is computed at time $T = 0.2$ with 100 cells and a reference solution with the second order scheme using 32000 cells. The functions $\alpha$ and $\beta$ are chosen according to the isothermal atmosphere from Section 2.1. In Figure 1, the pressure perturbation $p(x) - p_0(x)$ and the velocity perturbation are plotted in comparison with the initial perturbation. One can observe, that the second order scheme captures the peaks of the resulting waves more accurately than the first order scheme.

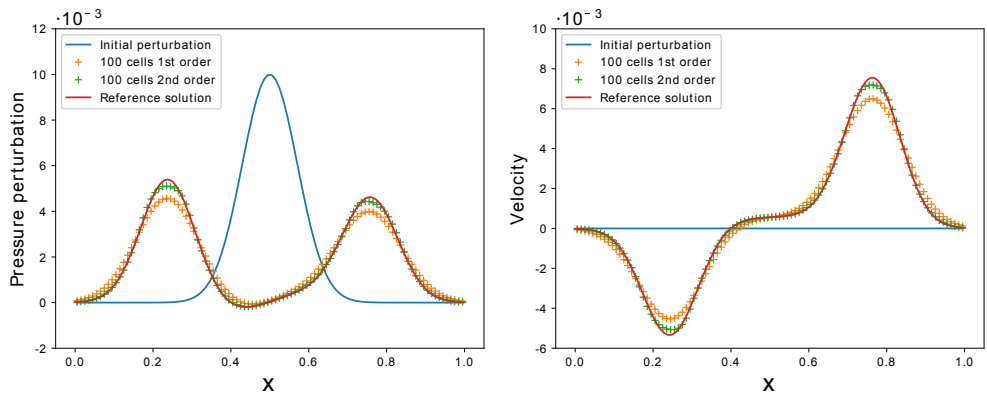

Figure 1: Perturbation in pressure (left) and in velocity (right).

### 6.4 Rayleigh-Taylor instability

In this test case taken from,[9] a perturbation in density is put over an isothermal solution in polar coordinates with potential $\Phi = r$. The initial pressure and density on the domain $[-1, 1] \times [-1, 1]$ are given by

$$p = \begin{cases} \exp(-r) & r \leq r_0 \\ \exp(-\frac{r}{\mu} + r_0 \frac{1-\mu}{\mu}) & r > r_0 \end{cases}, \quad \rho = \begin{cases} \exp(-r) & r \leq r_i(\theta) \\ \frac{1}{\mu} \exp(-\frac{r}{\mu} + r_0 \frac{1-\mu}{\mu}) & r > r_i(\theta) \end{cases}, \quad (6.6)$$

where $r_i(\theta) = r_0(1 + \nu \cos(k\theta))$ and $\mu = \exp(-r_0)/\exp(-r_0) + \Delta_\rho)$. This results in a jump in density by an amount of $\Delta_\rho$ at the interface defined by $r = r_i$ whereas the pressure is continuous. Following,[16] we take $\Delta_\rho = 0.1$, $\eta = 0.02$, $k = 20$. For computation, we use a mesh of $240 \times 240$ cells. For $r < r_0(1 - \eta)$ and $r > r_0(1 + \eta)$, the initial condition is in stable equilibrium but due to the discontinuous density, a Rayleigh-Taylor instability develops at the interface defined by $r = r_i$. In Figure 2 the initial value and the solution at times $t = 2.9, 3.8$ and $t = 5.0$ for the density are plotted. It can be seen that instabilities occur only around the discontinuous interface since a well-balanced scheme is used for the computations. Since a Cartesian mesh and dimensional splitting is used in the computations, the instabilities evolve along the coordinate axes and thus the solution is not radial.

### 6.5 Rarefraction test

In order to demonstrate the positivity preserving property, we follow the 1-2-3 rarefraction test in[10] where we set $\rho$ and $p$ isothermal with a quadratic potential $\Phi(\mathbf{x}) = \frac{1}{2}[(x_1 - 0.5)^2 + (x_2 - 0.5)^2]$ centred around $\mathbf{x} = (0.5, 0.5)$. As initial velocity we set

$$u_1 = \begin{cases} -2 & \text{for } x_1 < 0.5, \\ 2 & \text{for } x_1 \geq 0.5, \end{cases} \quad u_2 = 0 \quad (6.7)$$

The calculations were performed with the second order scheme on 100 cells on the domain $[0, 1]^2$ up to $T_f = 0.1$. We want to mention, that throughout the computations the limiting procedure described in[14] to avoid negative values in the reconstruction was not necessary. As can be seen from Figure 3, the pressure, energy and density come close to 0 but remains positive throughout the simulation.

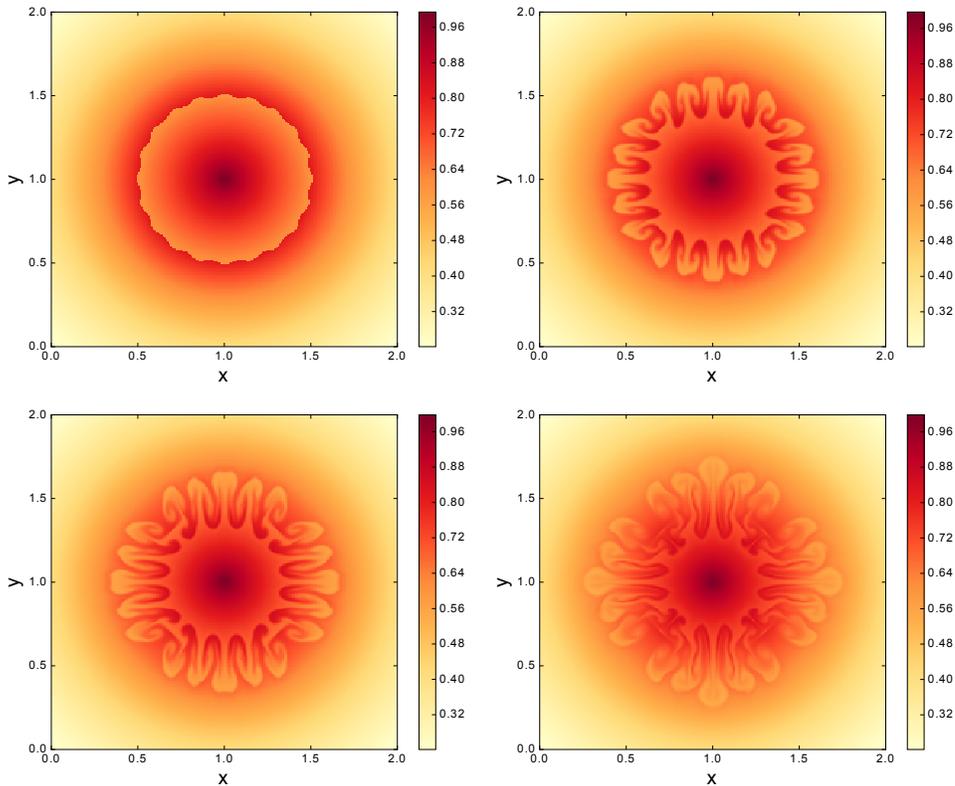

Figure 2: Rayleigh-Taylor instability in density in radial gravitational field at times $t = 0$ (top left), $t = 2.9$ (top right), $t = 3.8$ (bottom left) and $t = 5.0$ (bottom right).

# 7 Conclusion

A method for well-balancing arbitrary given hydrostatic equilibria of the compressible Euler equations with gravity was presented. It has been combined with a relaxation solver based on finite volume discretization. For this scheme, a second order well-balanced extension was described by using linear reconstruction and transformation of the pressure variables. The resulting scheme was given in a higher dimensional framework. To validate the proven properties, well-balancednes, robustness and second order accuracy, numerical examples were shown. Therein a novel higher dimensional exact solution (6.5) was used which allows to ascertain the numerical error well. An additional feature important for physical computations is the robustness of the presented scheme. It guarantees the positivity of density and total energy for the first as well as for the second order scheme.

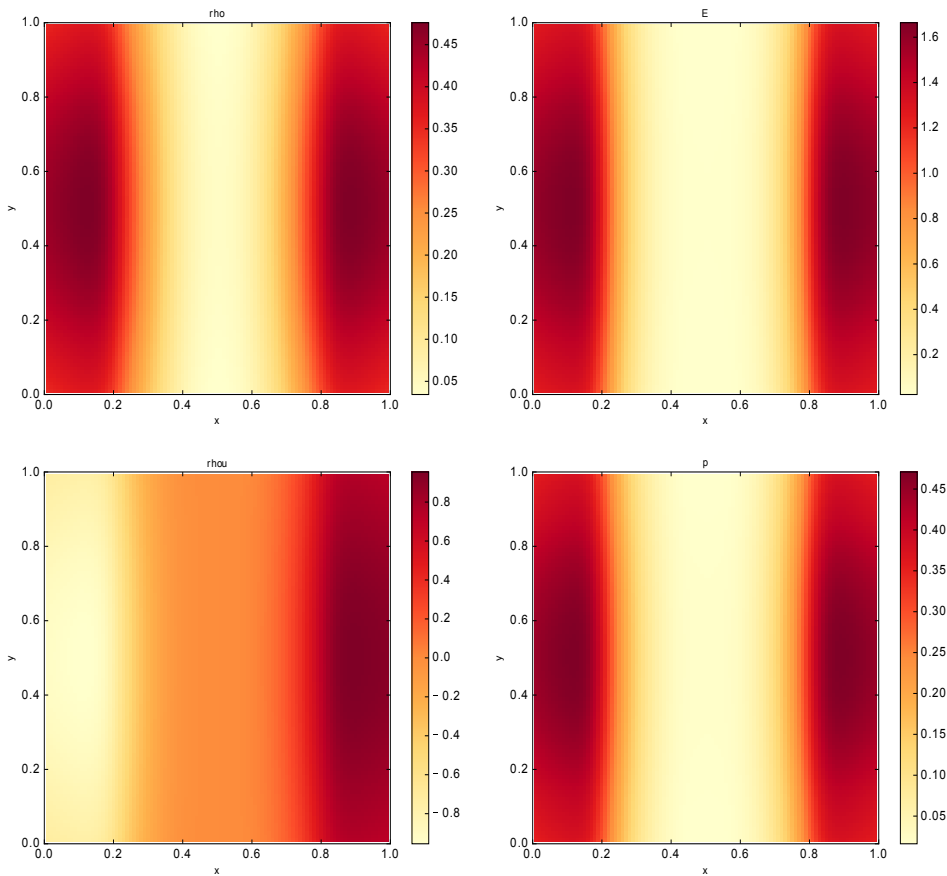

Figure 3: Density, energy, momentum and pressure for the rarefraction test with initial velocity $u = \pm 2.0$ at $T_f = 0.1$.

# Acknowledgements

The first author would like to thank INDAM-DP-COFUND-2015, grant 713485 for financial support.

|   | 50 | 100 | 200 | 400 | 800 | 1600 |
|---|---|---|---|---|---|---|
| $\rho$ | 9.083e-06 | 1.673e-06 | 3.428e-07 | 7.632e-08 | 1.790e-08 | 4.328e-09 |
|   | — | 2.441 | 2.287 | 2.167 | 2.092 | 2.048 |
| $\rho u$ | 5.171e-05 | 1.303e-05 | 3.268e-06 | 8.182e-07 | 2.047e-07 | 5.120e-08 |
|   | — | 1.989 | 1.995 | 1.998 | 1.999 | 1.999 |
| $E$ | 1.841e-05 | 3.433e-06 | 7.114e-07 | 1.595e-07 | 3.757e-08 | 9.106e-09 |
|   | — | 2.423 | 2.271 | 2.157 | 2.086 | 2.045 |

Table 1: $L^1$ error and convergence rates at $T_f = 1.0$ balancing (6.2) with an isothermal equilibrium using $\overline{S}$.

| N | $\rho$ | $\rho u$ | $\rho v$ | $\rho w$ | $E$ |
|---|---|---|---|---|---|
| 50 | 5.996E-017 | 1.450E-016 | 1.450E-016 | 1.438E-016 | 7.404E-017 |
| 150 | 4.930E-017 | 2.121E-016 | 2.121E-016 | 2.128E-016 | 6.426E-017 |

Table 2: $L^1$ error with respect to initial values in density, momentum and energy for isothermal equilibrium of Section 2.1 at $T_f = 1.0$.

| N | $\rho$ | $\rho u$ | $\rho v$ | $E$ |
|---|---|---|---|---|
| 100 | 4.796E-017 | 1.188E-016 | 1.188E-016 | 8.684E-017 |
| 500 | 1.036E-016 | 5.525E-016 | 5.525E-016 | 2.209E-016 |

Table 3: $L^1$ error with respect to initial values in density, momentum and energy for a polytropic equilibrium of section 2.2 at $T_f = 1.0$.

| N    | $\rho$     | $\rho u$   | $E$        |
|------|------------|------------|------------|
| 200  | 2.187E-016 | 3.742E-015 | 1.263E-015 |
| 1000 | 2.220E-018 | 5.149E-016 | 4.529E-017 |

| N   | $\rho$     | $\rho u$   | $\rho v$   | $E$        |
|-----|------------|------------|------------|------------|
| 100 | 9.294E-017 | 2.285E-015 | 2.285E-015 | 1.214E-015 |
| 500 | 5.571E-016 | 6.459E-015 | 6.459E-015 | 3.257E-015 |

| N   | $\rho$     | $\rho u$   | $\rho v$   | $\rho w$   | $E$        |
|-----|------------|------------|------------|------------|------------|
| 50  | 2.137E-016 | 3.092E-015 | 3.098E-015 | 3.032E-015 | 2.174E-015 |
| 150 | 4.930E-015 | 1.292E-014 | 1.299E-014 | 1.298E-014 | 5.000E-014 |

Table 4: $L^1$ error with respect to initial values in density, momentum and energy for a general stationary state at $T_f = 1.0$.

| N | $\rho$ | | $\rho u$ | | $\rho v$ | | $E$ | |
|---|---|---|---|---|---|---|---|---|
| 50  | 3.223E-004 | —     | 6.438E-003 | —     | 6.438E-003 | —     | 1.290E-001 | —     |
| 100 | 8.604E-005 | 1.905 | 1.716E-003 | 1.907 | 1.716E-003 | 1.907 | 3.444E-002 | 1.905 |
| 150 | 3.983E-005 | 1.899 | 7.936E-004 | 1.901 | 7.936E-004 | 1.901 | 1.594E-002 | 1.899 |
| 200 | 2.298E-005 | 1.910 | 4.573E-004 | 1.915 | 4.573E-004 | 1.915 | 9.206E-003 | 1.909 |
| 250 | 1.497E-005 | 1.919 | 2.975E-004 | 1.926 | 2.975E-004 | 1.926 | 6.000E-003 | 1.919 |
| 300 | 1.058E-005 | 1.904 | 2.099E-004 | 1.912 | 2.099E-004 | 1.912 | 4.244E-003 | 1.898 |
| 350 | 7.885E-006 | 1.909 | 1.562E-004 | 1.919 | 1.562E-004 | 1.919 | 3.166E-003 | 1.901 |
| 400 | 6.095E-006 | 1.928 | 1.205E-004 | 1.937 | 1.205E-004 | 1.937 | 2.453E-003 | 1.912 |

Table 5: $L^1$ error and convergence rates with respect to the exact solution in 2D, $T_f = 0.01$.

| N | $\rho$ | | $\rho u$ | | $\rho v$ | | $\rho w$ | | $E$ | |
|---|---|---|---|---|---|---|---|---|---|---|
| 25  | 5.209E-006 | —     | 2.031E-006 | —     | 2.031E-006 | —     | 2.031E-006 | —     | 1.215E-005 | —     |
| 50  | 1.739E-006 | 1.582 | 5.346E-007 | 1.926 | 5.346E-007 | 1.926 | 5.346E-007 | 1.926 | 4.205E-006 | 1.531 |
| 75  | 8.507E-007 | 1.763 | 2.379E-007 | 1.996 | 2.379E-007 | 1.996 | 2.379E-007 | 1.996 | 2.077E-006 | 1.739 |
| 100 | 5.013E-007 | 1.838 | 1.332E-007 | 2.014 | 1.332E-007 | 2.014 | 1.332E-007 | 2.014 | 1.228E-006 | 1.824 |
| 125 | 3.295E-007 | 1.879 | 8.478E-008 | 2.026 | 8.478E-008 | 2.026 | 8.478E-008 | 2.026 | 8.093E-007 | 1.871 |
| 150 | 2.330E-007 | 1.900 | 5.871E-008 | 2.015 | 5.871E-008 | 2.015 | 5.871E-008 | 2.015 | 5.730E-007 | 1.893 |
| 175 | 1.732E-007 | 1.923 | 4.291E-008 | 2.033 | 4.291E-008 | 2.033 | 4.291E-008 | 2.033 | 4.262E-007 | 1.919 |
| 200 | 1.339E-007 | 1.929 | 3.280E-008 | 2.015 | 3.280E-008 | 2.015 | 3.280E-008 | 2.015 | 3.297E-007 | 1.925 |

Table 6: $L^1$ error and convergence rates with respect to the exact solution in 3D, $T_f = 0.01$.